\documentclass[11pt]{article}
\usepackage{epsf,amssymb,latexsym,amsmath}
\usepackage{graphicx}
\usepackage{tikz}
\usepackage[autostyle=true]{csquotes}

\newcommand{\proof}{\noindent{\bf Proof.\ }}
\newcommand{\qed}{\hfill $\square$ \bigskip}

\newtheorem{theorem}{\bf Theorem}[section]
\newtheorem{corollary}[theorem]{\bf Corollary}
\newtheorem{lemma}[theorem]{\bf Lemma}
\newtheorem{proposition}[theorem]{\bf Proposition}

\newtheorem{definition}[theorem]{\bf Definition}

\newtheorem{problem}[theorem]{\bf Problem}

\begin{document}

\title{Polyhedra without cubic vertices are prism-hamiltonian}

\author{
Simon \v Spacapan\footnote{ University of Maribor, FME, Smetanova 17,
2000 Maribor, Slovenia. e-mail: simon.spacapan @um.si.
}}
\date{\today}

\maketitle

\begin{abstract} 
 
The prism over a graph $G$ is the Cartesian product of $G$ with the complete graph on two vertices. 
A graph $G$ is prism-hamiltonian if the prism over $G$ is hamiltonian. 
We prove that every polyhedral graph (i.e. 3-connected planar graph) of minimum degree at least four is prism-hamiltonian. 
\end{abstract}



\section{Introduction}

The study of hamiltonicity of planar graphs is largely concerned with finding subclasses of 3-connected planar graphs for which each member of the subclass is hamiltonian or has some hamiltonian-type property. 
One such result was obtained in 1956  by Tutte who proved that all $4$-connected planar graphs are hamiltonian \cite{tutte2}. Although not every 3-connected planar graph is hamiltonian it is possible to prove that this class of graphs satisfies (hamiltonian-type) properties weaker than hamiltonicity.  A  2-walk in a graph is a closed spanning walk that visits every vertex at most twice.  
Clearly, every hamiltonian graph has a 2-walk. In  \cite{richter} Gao and Richter proved that every $3$-connected planar graph has a 2-walk.

There is an extensive list of non-hamiltonian 3-connected planar graphs with special properties, such as 
graphs with small order and size   \cite{bar}, 
  plane triangulations  \cite{zamfirescu}, regular graphs 
 \cite{tutte1, zamfirescu1},  $K_{2,6}$-minor-free  graphs  \cite{eli1}, 
 and graphs with few 3-cuts  \cite{brink2}.  
However  some   classes of graphs mentioned above are prism-hamiltonian. For example every plane triangulation is prism-hamiltonian \cite{bib}, and every cubic 3-connected graph is prism-hamiltonian \cite{kaiser},\cite{paul}. It is well known that every prism-hamiltonian graph has a 2-walk, 
so the result obtained in \cite{bib} strengthens the result of Gao and Richter mentioned above.

Rosenfeld and Barnette  \cite{domneva}  conjectured    that every 3-connected planar graph is prism-hamiltonian (see also  \cite{kral}).
This conjecture was recently refuted in \cite{jaz} where vertex degrees play a central role in construction of counterexamples. In particular every counterexample to  Rosenfeld-Barnette conjecture 
given in \cite{jaz} has many cubic vertices and two vertices of \enquote{high} degree (linear in order of the graph). 
In \cite{zam} the authors show that there is an infinite family of 3-connected planar graphs, each of them not prism-hamiltonian, such that the ratio of cubic vertices tends to 1 when the order goes to infinity, and  
maximum degree stays bounded by 36. 

Vertex degrees in relation to hamiltonicity properties are discussed already by Ore in \cite{ore} and later by 
Jackson and Wormald in \cite{jackson}. Let $\sigma_k(G)$ be the minimum sum of vertex degrees of an independent set of $k$ vertices. Ore showed that $\sigma_2(G)\geq n$ implies  that $G$ is hamiltonian, and 
Jackson and Wormald showed that $\sigma_3(G)\geq n$ implies that  $G$ has a 2-walk  (provided that $G$ is connected). This was strenghtened by Ozeki in \cite{kenta} who showed that  $\sigma_3(G)\geq n$  implies that  $G$ is prism-hamiltonian.


In this paper we prove that every 3-connected planar graph of minimum degree at least four is prism-hamiltonian.  
 Equivalently, every    3-connected planar graph which is not prism-hamiltonian must have  
at least one cubic vertex.  In particular this implies that every regular 3-connected planar graph is prism-hamiltonian. 
The class of 3-connected planar graphs of minimum degree at least four is neither 
hamiltonian nor traceable (even when restricted to plane triangulations, or to regular graphs), see \cite{zamfirescu} and \cite{zamfirescu1}. In this sense prism-hamiltonicity appears to be the strongest hamiltonian-type property this class has. 

The proof  we give in this article builds on results obtained in \cite{richter}, where a 
 method of decomposing graphs into plain chains is developed. In  \cite{richter} the authors work with circuit graphs (which where originally defined in \cite{barnette}). A plane graph is a circuit graph if it is obtained from a 3-connected plane graph $G$ by deleting all vertices that lie in the exterior of a cycle of $G$. A cactus is a connected graph $G$ such that every block of $G$  is either a $K_2$ or a cycle, and such that  every vertex of $G$ is contained in at most two blocks of $G$ (the last condition is usually omitted, however for us it will be crucial, so we include it in the definition).
The main result of 
\cite{richter} is that any circuit graph 
(and hence also any 3-connected plane graph) has a spanning 
cactus as a subgraph. 
Here we improve this result by proving that any circuit graph   
with no internal cubic vertex  has a spanning bipartite cactus as a subgraph. Every 
cactus has a 2-walk while every bipartite cactus is prism-hamiltonian. Our result thus implies that 
 circuit graphs with all internal 
vertices of degree at least 4 are prism-hamiltonian.

We mention that 3-connected planar graphs of minimum degree at least 4 
also appear in \cite{thomassen} where the author proved that 
no graph in this class is hypohamiltonian.

\section{Preliminaries}
We refer to \cite{mt} for terminology not defined here. 
Let $G=(V(G),E(G))$ be a graph, $x\in V(G)$ and $X\subseteq V(G)$. 
We say that $x$ is {\em adjacent} to $X$, 
if $x$ is adjacent to some vertex of $X$. 
If $u$ and $v$ are adjacent  then $e=uv$ denotes the edge with endvertices $u$ and $v$; the subgraph induced by $u$ and $v$ is a path  denoted by $u,v$. 
The {\em union} of graphs $G=(V(G),E(G))$ and  $H=(V(H),E(H))$ is the graph 
$G\cup H=(V(G)\cup V(H),E(G)\cup E(H))$ and the {\em intersection} of $G$ and $H$ is  
$G\cap H=(V(G)\cap V(H),E(G)\cap E(H))$. The graph $G-X$ is  obtained from $G$ 
by deleting all vertices in $X$ and   edges incident to a vertex in $X$. Similarly, for 
$M\subseteq E(G)$, $G-M$ is the graph obtained from $G$ 
by deleting all edges in $M$. If $X=\{x\}$ we write $G-x$ instead of $G-\{x\}$.

Let $G$ be a plane graph.  
 Vertices and edges  incident to the unbounded face of $G$ are called  {\em external vertices} and {\em  external edges}, respectively.  
If a vertex (or an edge) is not an external vertex (or edge), then 
 it is called an  {\em internal vertex} (or an {\em  internal edge}). 
A path $P$ is an {\em external} resp. {\em internal}  path of $G$ if all edges of $P$ are external resp. internal edges.

We use $[n]$ to denote the set of positive integers less or equal $n$. 
A path of odd/even length is called an {\em odd/even path}, respectively.  
Similarly we define {\em odd} and {\em even faces}, based on the parity of their degree.

Recall that every vertex of a cactus $G$ is contained in 
at most two blocks of $G$. A vertex of a catus $G$ is {\em good} if it is contained in exactly one block of $G$.  

A {\em prism} over a graph $G$ is the Cartesian product of $G$ and the complete graph on two vertices $K_2$. The following proposition is given in \cite{eli3} (Theorem 2.3.). 
For the sake of completness we include the proof of it also here. 

\begin{proposition} \label{bolje}
Every bipartite  cactus is prism-hamiltonian.   
\end{proposition}
\proof 
We denote  $V(K_2)=\{a,b\}$.  We use induction to prove the following stronger statement. Every prism $G\Box K_2$ over a bipartite cactus $G$ has a Hamilton cycle $C$ such that for  
every good vertex $x$ of $G$, we have $(x,a)(x,b)\in E(C)$.
This is clearly true when $G$ is an even  cycle or $K_2$. 

Let $G$ be a bipartite cactus and assume that the statement is true for all   bipartite cactuses with fewer vertices than 
$|V(G)|$. If all vertices of $G$ are good, then $G$ is an even cycle or $K_2$. Otherwise, there  is a  vertex $u$, which is not a good vertex of $G$. 
Hence, $u$ is contained in exactly two blocks of $G$.

Let $G_1'$ and $G_2'$ be connected components of $G-x$, and let $G_1=G-G_2'$ and $G_2=G-G_1'$.
Both, $G_1$ and $G_2$, are   bipartite cactuses. Moreover, $x$ is a good vertex in $G_i$, for $i=1,2$. 
By  induction hypothesis there is a Hamilton cycle $C_i$ in $G_i$ such that $C_i$ uses the edge $e=(x,a)(x,b)$ in $G_i$. 
The desired Hamilton cycle in $G$ is $(C_1\cup C_2)-e$. Observe that every good vertex of $G$ is a good vertex of $G_1$ or $G_2$. 
It follows that for  
every good vertex $x$ of $G$, we have $(x,a)(x,b)\in E(C)$.
\qed

\begin{corollary}
Every graph $G$, that has  a bipartite cactus  $H$ as a spanning subgraph,  is prism-hamiltonian. 
\end{corollary}

If  $G$ is plane graph and $H$ is a subgraph of $G$, then $H$ is also a plane graph and we assume that the embedding of $H$ in the plane is the one given by $G$.  

Let $G$ be a plane graph and  $G^+$ the graph obtained from 
$G$ by adding a vertex to $G$ and making it adjacent to all external vertices of $G$. 
The graph $G$ is a {\em circuit graph} if  $G^+$ is 3-connected. 

It follows from the definition that any circuit graph is 2-connected, and hence every face of a circuit graph is bounded by a cycle. 
If $G$ is a 3-connected plane graph  (or if $G$ is a circuit graph) and $C$ is a cycle of $G$, then the subgraph of $G$ bounded by $C$ is a circuit graph. 
Observe also that for any circuit graph $G$ with outer cycle $C$, and any separating set $S$ of size 2 in $G$, every connected component of $G-S$ 
intersects $C$.

A graph $G$ is a {\em chain of blocks} if the block-cutvertex graph of $G$ is a path. We denote the blocks and cutvertices of $G$, by 
$$B_1,b_1,B_2,\ldots,b_{n-1}, B_n\,,$$ 
where $B_i$  are blocks for $i\in [n]$, and $b_i\in V(B_i)\cap V(B_{i+1})$ are cutvertices of $G$ for $i\in [n-1]$. 
A plane graph $G$ is a {\em plane chain of blocks} if it is a 
chain of blocks $$G=B_1,b_1,B_2,\ldots,b_{n-1}, B_n$$ such that every external vertex of $B_i, i \in [n]$ is  also an external vertex of $G$. The following lemma is given in \cite{richter} (Lemma 3, p. 261). 

\begin{lemma}\label{plainchain}
Let $G$ be a circuit graph with outer cycle $C$ and let $x\in V(C)$. 
Let $x'$ and $x''$ be the neighbors of $x$ in $C$. Then 
\begin{itemize}
\item[(i)] $G-x$ is a plane chain of blocks $B_1,b_1,B_2,\ldots,b_{n-1}, B_n$ and each nontrivial block of $G-x$ is a circuit graph. 
\item[(ii)]Setting  $x'=b_0$ and $x''=b_n$, then $ B_i\cap C$ is a path in $C$ with  endvertices $b_{i-1}$ and $b_i$,  for every $i\in [n]$.  
 \end{itemize} 
\end{lemma}

It follows from the above lemma that for every  nontrivial block $B_i$ of $G-x$, with outer cycle $C_i$, 
 $C_i$ is the union of two $b_{i-1}b_i$-paths $P_i$ and $P_i'$, where 
$P_i$ is an internal path in $G$ and $P_i'$ is an external path in $G$.

\section{The proof of main result}

In this section we prove that any circuit graph $G$ such that every internal vertex of 
$G$ is of degree at least 4 is prism-hamiltonian.

\begin{definition}
Let $G$ be a circuit graph with outer cycle $C$ and let $x,y\in V(C)$. We say that $G$ is  bad with respect to $x$ and $y$ if 
\begin{itemize}
\item[(i)] $G$ has exactly one bounded odd face $F$
\item[(ii)] $x$ and $y$ are incident to $F$
\item[(iii)] If $x$ and $y$ are adjacent, then $e=xy$ is an internal edge of $G$.
 \end{itemize} 
We say that $G$ is  good with respect to $x$ and $y$  if it's not  bad with respect to $x$ and $y$.  
\end{definition}

If $G$ is a circuit graph and $G$ is bad with respect to $x$ and $y$, then there is no hamiltonian cycle $C$ in 
$G\Box K_2$ such that $C$ uses   vertical edges at $x$ and $y$ (edges between the two layers of $G$). 
For example, an odd cycle is bad with respect to any two non-adjacent vertices, and hence the prism over 
an odd cycle has no hamiltonian cycle that uses vertical edges at  two non-adjacent vertices of this cycle.  
Conversely, it turns out (and is a consequence of Theorem \ref{glavni}) that for any circuit graph $G$ with all internal vertices of degree at least 4, and any external vertices $x$ and $y$ of $G$ such that 
$G$ is good with respect to $x$ and $y$, there is a hamiltonian cycle in 
$G\Box K_2$ that uses   vertical edges at $x$ and $y$. 

Note also that   a bipartite circuit graph $B$ is good with respect to any two external vertices of $B$ 
(this fact we shall use frequently). 
In order to simplify the formulation of statements, we 
 also say that complete graphs $K_1$ and $K_2$ are good with respect to any of its vertices.

\begin{definition}
Let $G=B_1,b_1,B_2,\ldots,b_{n-1}, B_n$ be a plain chain of blocks such that each nontrivial block $B_i$ is a circuit graph. 
Let $b_0\neq b_1$ be an external vertex of $B_1$, and $b_n\neq b_{n-1}$ be an external vertex of $B_n$. We say that 
$G$ is a good  chain with respect to $b_0$ and $b_n$ if  $B_i$ is good with respect to $b_{i-1}$ and $b_i$ 
for every  $i\in [n]$. 
\end{definition}

The same definition is used when only one of the two vertices $b_0$ and $b_n$ is given, 
and in this case we say that $G$ is a good chain with respect to $b_0$ or with respect to $b_n$.  
If $G=B_1$ has only one block   we say that $G$ is a 
good chain with respect to any external vertex of $G$. 
 
\begin{lemma} \label{dvodelen}
Let $B$ be a bipartite circuit graph with outer cycle $C$ such that all internal vertices of $B$ are of degree at least 4. 
Then $B$ has at least 4 external vertices of degree 2. 
\end{lemma}
\proof
Let $C$ be a $k$-cycle, $k\geq 4$. Let ${\mathcal F}$ be the set of faces of $B$, and set $e=|E(B)|, v=|V(B)|$ and $f=|{\mathcal F}|$. Since $B$ is bipartite   
$$2e=\sum_{F\in {\mathcal F}} \deg (F)\geq 4(f-1)+k\,.$$
We use the Euler's formula to obtain 
$$\sum_{x\in V(B)}\deg(x)=2e\leq 4v-k-4\,.$$
Since every internal vertex of $B$ is of degree at least 4 we get 
$$\sum_{x\in V(C)}\deg(x)\leq 3k-4\,.$$ 
Since all vertices of $C$ are of degree at least 2, the claim of the lemma follows from the pigeonhole principle. 
\qed

The following lemma is a well known fact. 

\begin{lemma}\label{sodalica}
A plane graph $G$ is bipartite if and only if all bounded faces of $G$ are even. 
\end{lemma}

\begin{lemma}\label{osnovna1}
Let $B$ be a circuit graph with outer cycle $C$ such that all internal vertices of $B$ are of degree at least 4. Let $x$ and $y$ be any vertices of $C$, and $Q$ a 
$xy$-path in $C$. Suppose that   
all vertices in  $V(C)\setminus V(Q)$ are of degree at least three in $B$.  
Then $B$ is good with respect to $x$ and $y$.  
\end{lemma}

\proof 
Suppose to the contrary, that $B$ is bad with respect to $x$ and $y$. 
Then  $x$ and $y$ are incident to odd face $F$ of $B$, and $F$ is the only bounded odd face of $B$. 
Moreover, $x$ and $y$ are not adjacent in $C$.  It follows that 
$B-\{x,y\}$ has exactly two components.

Let $H$ be the component of $B-\{x,y\}$ that contains a vertex of $Q$. If $xy\in E(B)$   and 
 $F$ is contained  in the exterior of the cycle   $E(Q)\cup \{xy\}$ define $H'=(B-V(H))-xy$. 
Otherwise  define $H'=B-V(H)$.   
$H'$ is a plain chain of blocks and each nontrivial block of $H'$ is a bipartite circuit graph 
(by Lemma \ref{sodalica}), so assume 
$$H'=D_1,d_1,D_2,\ldots,d_{m-1}, D_m\,.$$
Let $d_0=x$ and $d_m=y$.
If  $j\in [m]$ and  $u\in V(D_j)\setminus \{d_{j-1},d_j\}$, then 
$\deg_{D_j}(u)>2$. 
So if $D_j$ is nontrivial, then it has at most two vertices of degree 2 in $D_j$; since $D_j$ is bipartite this contradicts Lemma \ref{dvodelen}. It follows that all blocks of $H'$ 
are trivial. If $H'$ is $K_2$ then $x$ and $y$ are adjacent in $C$ 
(a contradiction), otherwise a vertex in $V(C)\setminus V(Q)$ is of degree  $\leq 2$ (this contradicts the assumption of the lemma).   
\qed

\begin{lemma}\label{osnovna}
Let $B$ be a circuit graph with outer cycle $C$ such that all internal vertices of $B$ are of degree at least 4. Let $x\in V(C)$ be any vertex 
and 
$$B-x=B_1,b_1,B_2,\ldots,b_{n-1}, B_n\,.$$ 
Let $b_0\in V(B_1)$ and $b_n\in V(B_n)$  be the neighbors of $x$ in $C$. Then
for every $i\in [n]$, $B_i$ is good with respect to $b_{i-1}$ and $b_i$. 

\end{lemma}

\proof
Let $B_i$ be a nontrivial block with outer cycle $C_i$, and define $Q=C\cap B_i$. 
$Q$ is a path in $C_i$ with endvertices $b_{i-1}$ and $b_i$, and 
every vertex in  $V(C_i)\setminus V(Q)$ is of degree more than 2 in $B_i$. By Lemma 
\ref{osnovna1}, $B_i$ is good with respect to  $b_{i-1}$ and $b_i$.   
\qed

\begin{lemma}\label{osnovna2}
Let $B$ be a  circuit graph with outer cycle $C$ such that all internal vertices of $B$ are of degree at least 4. Let $x$ and $y$ be any vertices of $C$  and $Q$ a  
$xy$-path in $C$ such that   
all vertices  in  $V(C)\setminus V(Q)$ are of degree at least three in $B$.  If $B-x$ is bipartite, then 
$|V(B_i\cap Q)|\geq 2$ for every block $B_i$ of $B-x$.   
\end{lemma}

\proof
If $V(C)=V(Q)$, the lemma follows from Lemma \ref{plainchain}. Assume $V(C)\neq V(Q)$, and 
let $u\in V(C)\setminus V(Q)$ be the neighbor of $x$. 
The block $B_1$ of $B-x$ containing $u$ is nontrivial, for otherwise $\deg_B(u)=2$. 
If $|V(B_1\cap Q)|< 2$, then $B_1$ has at most two vertices of degree two in $B_1$. Therefore, by Lemma \ref{dvodelen}, $B_1$ is non-bipartite and hence $B-x$ is non-bipartite.
\qed

\begin{lemma} \label{skupek}
Let $B$ be a circuit graph with outer cycle $C$ such that all internal vertices of $B$ are of degree at least 4. Let $x\in V(C)$ be any vertex, 
and let
$$B-x=B_1,b_1,B_2,\ldots,b_{n-1}, B_n\,.$$
Then for every $k \in [n-1]$ the graph
$$G=B-\bigcup_{i=k+1}^n V(B_i)$$
is a good chain with respect to $x$. 

\end{lemma}

\proof
Let $b_0\in V(B_1)$ be the neighbor of $x$ in $C$. Denote  the path $x, b_0$ by $B_0$.  

{\em Case 1: 
 Suppose that $B_k$ is trivial. } \\
 If $x$ is not adjacent to a vertex in $G-b_0$, then  
 $G$ induces a plain chain of blocks
$$B_0,b_0,B_1,\ldots,b_{k-2}, B_{k-1}\,$$
and, by Lemma \ref{osnovna}, $B_i$ is good with respect to $b_{i-1}$ and $b_i$ for $i\in [k-1]$. 
Assume therefore that $x$ is adjacent to a vertex in $G-b_0$.  Let $\ell\in [k]$ be the maximum number such that $x$ is  adjacent to $B_\ell-\{b_{\ell-1},b_k\}$.  
Since $B_k$ is trivial, $\ell\neq k$.  
The graph $H$ induced by $\bigcup_{i=0}^\ell V(B_i)$
is a nontrivial block of $G$. Moreover, since $H$ is a subgraph of $B$ bounded by a cycle of $B$, $H$ is a circuit graph.  
Note also that if  $x$ and $b_\ell$ are incident to a bounded face 
$F$ of $H$, then $x$ and $b_\ell$ are adjacent, moreover  
$xb_\ell$ is an external edge of $H$.  It follows that 
  $H$ is good with respect to $x$ and $b_\ell$, and therefore 
$$G=H,b_\ell,B_{\ell+1},\ldots, b_{k-2}, B_{k-1}$$ 
is a good chain with respect to $x$.

{\em Case 2:  Suppose that $B_k$ is nontrivial. } \\
By Lemma \ref{plainchain}, $B_k-b_k$ is a plain chain of blocks, so let 
$$B_k-b_k=D_1,d_1,D_2,\ldots,d_{m-1}, D_m\,.$$
Let $C_k$ be the outer cycle of $B_k$, and $d_0\in V(D_1), d_m\in V(D_m)$  be the neighbors of $b_k$ in $C_k$. Without loss of generality assume that 
  $b_kd_0$ is an internal edge of $B$.  Note that $D_1$ is nontrivial if $b_{k-1}\neq d_0$, for otherwise 
$\deg_{B} (d_0)\leq 3$ (this is a contradiction because $d_0$ is an internal vertex of $B$ if $b_{k-1}\neq d_0$).

Suppose that $x$ is adjacent to   $D_1-\{b_{k-1},d_1\}$ (it is possible that $b_{k-1}=d_1$). Then 
$D_1$ is nontrivial. 
Let $j\in [m]$ be   such that $b_{k-1}\in V(D_j)\setminus \{d_{j-1}\}$ and let $H'$ be the graph induced by 
$$\bigcup_{i=0}^{k-1} V(B_i)\cup \bigcup_{i=1}^jV(D_i)\,.$$ 
$H'$ is bounded by a cylce of $B$, so it is a circuit graph. We shall prove that $H'$ is good with respect to $x$ and $d_j$. 
Suppose that $x$ and $d_j$ are incident to a face $F'$ of $H'$, and  that $F'$ is the only bounded odd face of $H'$. Then $d_j=b_{k-1}$, and 
all bounded faces of $D_1$ are even. This contradicts Lemma \ref{dvodelen}, because 
$\deg_{D_1}(u)>2$ for every $u\in V(D_1)\setminus \{d_0, d_1\}$. 
It follows that
 $$G=H',d_{j},D_{j+1},\ldots,d_{m-1},D_m\,$$
is a good chain with respect to $x$. 

Suppose that $x$ is not adjacent to $D_1-\{b_{k-1},d_1\}$. 
We claim that  $|V(D_1)\cap V(C)|\geq 2$.  
To prove the claim suppose the contrary, that $|V(D_1)\cap V(C)|<2$. Then $\{b_k,d_1\}$ is a separating set in $B$, and $D_1-\{b_k,d_1\}$ is a component 
of $B-\{b_k,d_1\}$ disjoint with $C$. 
It follows that  $B$ is not a circuit graph, a contradiction. This proves the claim.

Define $\ell$ and $H$ as in Case 1. 
We claim that   
 $$G=H,b_\ell,B_{\ell+1},\ldots, B_{k-1},b_{k-1},D_1,d_1,\ldots,d_{m-1},D_m$$ 
is a good chain with respect to $x$. 
We have already shown (in Case 1) that $H$ is good with respect to $x$ and $b_{\ell}$. 
By Lemma \ref{osnovna}, $B_i$ is good with respect to $b_{i-1}$ and $b_i$ for $i\in [k-1]\setminus [\ell]$, 
and $D_i$ is good with respect to $d_{i-1}$ and $d_i$ for $i\in [m], i\neq 1$. It remains to prove that 
$D_1$ is good with respect to $b_{k-1}$ and $d_1$. 
 Let $C'$ be the outer cycle of $D_1$ and let $Q=C\cap D_1$ (or equivalently  $Q=C\cap C'$).  
Note that for every vertex $z\in V(C')\setminus V(Q)$,   
$\deg_{D_1}(z)>2$. By Lemma \ref{osnovna1}, $D_1$ is good with respect to $b_{k-1}$ and $d_1$.
\qed

\begin{definition}\label{def}
Let $B$ be a circuit graph with outer cycle $C$. Let   $\{x,y\}\subseteq V(C)$ and $\{u_1,u_2\}\subseteq V(C)$ be any sets. 
A set of pairwise disjoint chains ${\mathcal C}=\{G_1,\ldots,G_k\}$  is a $(x,y;u_1,u_2)$-set of chains in $B$ if there exists a $xy$-path $P$ in $B$ such that 
\begin{itemize}
\item[(i)] $V(B)\setminus V(P)\subseteq \bigcup_{i=1}^{k}V(G_i)$ 
\item[(ii)] For $i\in [k]$, $G_i$ intersects $P$ in exactly one vertex $x_i$, and 
$G_i$ is a good chain with respect to $x_i$. 
\item[(iii)]  For $j\in [2]$, either $G_i$ is a good chain with respect to $u_j$ and $x_i$ for some $i\in [k]$, or $u_j \notin \bigcup_{i=1}^{k}V(G_i)$.  
\end{itemize}
 A  path $P$ that fulfills (i),(ii) and (iii)  is called a  ${\mathcal C}$-path. 
The set ${\mathcal C}$ is an odd or an even $(x,y;u_1,u_2)$-set of chains if there exits an 
odd or an even ${\mathcal C}$-path, respectively. 
\end{definition}

We say that a set of pairwise disjoint chains $G_1,\ldots,G_k$  is a $(x,y;u_1)$-set of chains
 if it satisfies (i),(ii) and (iii) for $j=1$. 
 Moreover, ${\mathcal C}$ is a $(x,y)$-set of chains if it satisfies  (i) and (ii)  of Definition \ref{def}. 

We also use Definition \ref{def} in slightly more general settings in which $B$ is a plain chain of blocks (and each block is a circuit graph). 
More precisely, if $B$ is a plain chain of blocks, $x,y$ are two external vertices of $B$, and ${\mathcal C}$ is a set of pairwise disjoint plain chains that 
satisfy (i),(ii) and (iii), then  ${\mathcal C}$ is a $(x,y;u_1,u_2)$-set of chains in $B$.


\begin{lemma} \label{spajanje}
 Let   $G$ be a bipartite plain chain of blocks 
$$G=B_1,b_1,B_2,\ldots,b_{n-1}, B_n $$ 
such that  for $i\in [n]$ each nontrivial block $B_i$ of $G$ is a circuit graph 
with outer cycle $C_i$. Suppose that $u,x,y\in V(C_j),u\neq b_j$, and that  
 ${\mathcal C}$ is a $(x,y;u,b_j)$-set of chains 
in $B_j$ for some $j\in [n]$. 
Then for every $\ell>j$ and any  vertex $v\in V(C_\ell)\setminus V(C_{\ell-1})$, there is a 
$(x,y;u,v)$-set of chains 
in $\bigcup_{i=j}^{\ell} B_i$. Moreover, if $u=b_{j-1}$, then 
for every  $\ell'<j$ and any vertex $v'\in V(C_{\ell'})\setminus V(C_{\ell'+1})$, there is a 
$(x,y;v',v)$-set of chains 
in $\bigcup_{i=\ell'}^{\ell} B_i$.
\end{lemma}

\proof
Let  $\ell>j$  and $v\in V(C_\ell)\setminus V(C_{\ell-1})$. Suppose that 
${\mathcal C}=\{ G_1,\ldots,G_k\}$ is a 
$(x,y;u,b_j)$-set of chains in $B_j$ and that $P$ is a ${\mathcal C}$-path. Then (a) or (b) occures.

\begin{itemize}
\item[(a)] There is a chain $G_{r} \in {\mathcal C}$ such that $G_{r}$ is a good chain with respect to $x_{r}$ and $b_{j}$, where $\{x_{r}\}=V(G_{r})\cap V(P)$
\item[(b)] $b_j \notin \bigcup_{i=1}^{k}V(G_i)$. 
\end{itemize}
In case (a), $G_{r}'=G_{r}\cup \bigcup_{i=j+1}^{\ell} B_i$ is a good chain with respect to $x_{r}$ and $v$, and therefore ${\mathcal C'} ={\mathcal C}\cup \{G_r'\}\setminus \{G_r\}$  
is a  $(x,y;u,v)$-set of chains in $ \bigcup_{i=j}^{\ell} B_i$. 
 In case (b),  $G_0=\bigcup_{i=j+1}^{\ell} B_i$ is a good chain with respect to $b_j$ and $v$, and therefore 
${\mathcal C'}=\{G_0,\ldots,G_{k}\}$ is a $(x,y;u,v)$-set of chains  in $ \bigcup_{i=j}^{\ell} B_i$.  In both cases 
a ${\mathcal C'}$-path is $P$. The last sentence of the lemma is proved analogously. 
\qed


If we use the notation of  the above lemma, we   note that a 
$(x,y;b_j)$-set of chains in $B_j$ can be extended to a 
$(x,y)$-set of chains in $\bigcup_{i=j}^{\ell} B_i$ (in fact the construction given in the above proof works also in this case). Note also that a $(x,y;u,v)$-set of chains in $\bigcup_{i=j}^{\ell} B_i$  exists also under the assumption that $B_i$ is bipartite for $i>j$ (and  $G$ may be non-bipartite). 
In \cite{richter} the following result was proved (Theorem 5, p.262).

\begin{theorem}\label{rihta}
Let $B$ be a bipartite circuit graph with outer cycle $C$. 
If $x,y\in V(C)$, then for any vertex $u\in V(C)$ (not necessarily distinct from $x$ and $y$) 
there exists  a $(x,y;u)$-set of chains in $B$. 
\end{theorem}

\begin{lemma}\label{posebna}
Let $B$ be a bipartite circuit graph with outer cycle $C$. Suppose that $x,y\in V(C)$ and that $Q$ is a $xy$-path in $C$. 
If every internal vertex of $B$ is of degree at least 4 and every vertex in $V(C)\setminus V(Q)$ is of degree at least 3 in $B$, then there exists  a
$(x,y;x,y)$-set of chains in $B$.  

\end{lemma}

\proof
By Lemma \ref{plainchain}, $B-x$ is a plain chain of blocks   
$$B-x=B_1,b_1,B_2,\ldots,b_{n-1}, B_n\,.$$
By Lemma \ref{osnovna2}, $|V(B_i\cap Q)|\geq 2$ for $i\in [n]$. We may assume, 
without loss of generality, that $y\in V(B_1)$ and $y\neq b_1$.
If $B_1$ is nontrivial then $B_1-y$ is a plain chain of blocks 
$$B_1-y=D_1,d_1,D_2,\ldots,d_{m-1}, D_m\,.$$ 
Let  $d_0\in V(D_1)$ be the neighbor of $y$ in $Q$, and define 
$k=\max\{i\,|\,D_i\cap Q\neq \emptyset\}$. 

{\em Case 1:} Suppose that $D_k$ intersects $Q$ in exactly one vertex (in this case $d_{k-1}$). Then $G=\bigcup_{i=k}^m D_i$ is a good chain with respect to $d_{k-1}$.

If $D_i$ is trivial define $P_i=D_i$ and ${\mathcal C_i}=\emptyset$, for $i\in [k-1]$.
If $D_i$ is nontrivial then, by Theorem \ref{rihta}, there is a $(d_{i-1},d_i;d_i)$-set of chains ${\mathcal C_i}$ in $D_i$, 
for $i\in [k-1]$. In this case let $P_i$ be a ${\mathcal C_i}$-path in $D_i$.

 If $B_i$ is trivial define $R_i=B_i$ and ${\mathcal F_i}=\emptyset$, for $i\in [n],i\neq 1$.
If $B_i$ is nontrivial then,  by Theorem \ref{rihta}, there is a $(b_{i-1},b_i;b_{i-1})$-set of chains ${\mathcal F_i}$ in $B_i$, 
 for $i\in[n],i\neq 1$.  In this case let $R_i$ be a ${\mathcal F_i}$-path in $B_i$.  Let $b_n$ be the  neighbor of $x$ in $Q$, 
and let $R_{n+1}$ be the path $x,b_n$. Additionally  let  $P_0$ be the path $y,d_0$.  
Define $$P=\bigcup_{i=0}^{k-1} P_i\cup  \bigcup_{i=2}^{n+1} R_i\,.$$
The chain $G$ together with chains ${\mathcal C_i}, i\in [k-1]$ and ${\mathcal F_i},i\in [n],i\neq 1$   is a 
$(x,y;x,y)$-set of chains  in $B$. If we call this set of chains ${\mathcal C}$, then $P$ is a ${\mathcal C}$-path.

{\em Case 2:} Suppose that $D_k$ intersects $Q$ in more than one vertex. Then, by   Theorem \ref{rihta}, there is a 
$(d_{k-1},b_1;d_k)$-set of chains ${\mathcal H}$ in $D_k$. By Lemma \ref{spajanje} (see also the note directly after Lemma \ref{spajanje}) there is a $(d_{k-1},b_1)$-set of chains in 
$\bigcup_{i=k}^m D_i$.  
The rest of the proof is similar as in Case 1. 

If $B_1$ is trivial, then $x$ and $y$ are adjacent in $C$ (for otherwise $\deg_B(u)=2$, where $u$ 
is the neighbor of $x$ in $V(C)\setminus V(Q)$) and $V(Q)=V(C)$. Define $R_1=B_1$. In this case 
$\bigcup_{i=2}^{n} {\mathcal F_i}$ is a  $(x,y;x,y)$-set of chains in $B$. The corresponding path is 
$\bigcup_{i=1}^{n+1} R_i.$
\qed

\begin{lemma}\label{to}
Let $B$ be a bipartite circuit graph with outer cycle $C$. Let $x,y,u_1,u_2\in V(C)$ be such that $\{x,y\}\neq \{u_1,u_2\}$. If every internal vertex of $B$ is of degree at least 4, then there is a 
$(x,y;u_1,u_2)$-set of chains in $B$. 
\end{lemma}

\proof
Suppose that the claim of the lemma is not true; let $B$ be a  counterexample with minimum number of vertices.   It's easy to verfy the lemma when $B$ is a 4-cycle, or any even cycle.

By Lemma \ref{plainchain}, $B-x$ is a plain chain of blocks
$$B-x=B_1,b_1,B_2,\ldots,b_{n-1}, B_n\,.$$
Let $Q$ and $Q'$ be the $xy$-paths in $C$.
 Let $b_0\in V(B_1)$ and $b_n\in V(B_n)$ be the neighbors of $x$ in $Q$ and $Q'$, respectively. We set  
$B_0=\emptyset$ (to avoid ambiguity in the following definitions).   
Let $k\in [n]$ be such that $y\in V(B_{k})\setminus V(B_{k-1})$, and let   
$k_j\in [n]$ be such that $u_j\in V(B_{k_j})\setminus  V(B_{k_j-1})$ for $j=1,2$ (if 
 $x\in  \{u_1,u_2\}$ this applies only to $k_1$ and we set $u_2=x$). We may assume, without loss of generality,  that $y\neq b_0$ (otherwise $y\neq b_n$ and we have a similar proof) and that $k_1\leq k_2$. 

We shall construct a $xy$-path $P$ in $B$. 
  For $i\in [k]$, if $B_i$ is trivial define ${\mathcal C_i}=\emptyset$ and $P_i=B_i$.
In the sequal we define ${\mathcal C_i}$ and $P_i$ for nontrivial blocks $B_i$. 

By minimality of $B$, Lemma  \ref{to} is true for every nontrivial block $B_i$ of $B-x$ 
and therefore, for every  $i\in [k]$ we can apply the statement of  Lemma \ref{to} to $B_i$. 

Denote the outer cycle of $B_i$ by $C_i$.  
Since $Q_i=B_i\cap C$ is a $b_{i-1}b_i$-path in $C_i$ and every vertex of $V(C_i)\setminus V(Q_i)$ is of degree at least 3 in $B_i$, we can also apply  Lemma \ref{posebna} to $B_i$. 
The following statements are obtained either by an application of Lemma   \ref{to} or Lemma \ref{posebna} to $B_i$. For $i\in [k-1]$ and $j=1,2$ there exists:

\begin{itemize}
\item[(i)] a $(b_{i-1},b_i;b_{i-1},b_i)$-set of chains in $B_i$ (by 
Lemma \ref{posebna}),
\item[(ii)]   a $(b_{k-1},y;b_{k-1},b_k)$-set of chains in $B_k$ (by minimality of $B$ (i.e. by the statement of Lemma \ref{to}) if $y\neq b_k$, and by 
Lemma \ref{posebna} if $y=b_k$),
\item[(iii)]    a $(b_{k_j-1},b_{k_j};b_{k_j-1},u_{j})$-set of chains in $B_{k_j}$ 
(by minimality of $B$ if $u_j\neq b_{k_j}$, and by Lemma \ref{posebna} if $u_j=b_{k_j}$),
\item[(iv)]   if $k_j=k$ and $u_j\neq y$, there is a $(b_{k-1},y;b_{k-1},u_j)$-set of chains in $B_k$ (by minimality of $B$),
\item[(v)]     if $k_j=k$  and $u_j\neq b_k$, there is a $(b_{k-1},y;u_j,b_{k})$-set of chains in $B_k$ (by minimality of $B$),
\item[(vi)]    if $k_1=k_2$,   there is a $(b_{{k_1}-1},b_{k_1};u_1,u_2)$-set of chains  in $B_{k_1}$ (by minimality of $B$),
\item[(vii)]   if $k_1=k_2=k$, there is a $(b_{k-1},y;u_1,u_2)$-set of chains in 
$B_k$ (by minimality of $B$). 
\end{itemize}

Since $\{x,y\}\neq \{u_1,u_2\}$ we may assume, without loss of generality, that  $y\notin \{u_1,u_2\}$. 
Therefore we have the following possibilities 
(1) $u_1,u_2\notin V(Q')$, (2) $u_1\notin V(Q'),u_2\notin V(Q)$,   (3) $u_2=x$ and $u_1\notin V(Q')$. 
All other possibilites are symmetric, and they can be obtained from  one of the above cases 
by exchanging the roles of $Q$ and $Q'$; for example, $u_1,u_2\notin V(Q')$ is symmetric to $u_1,u_2\notin 
V(Q)$.
Therefore we can also assume that $u_1\notin V(Q')$. With this assumption the following cases with regard to $k, k_1$ and $k_2$ may appear. Next to each particular case below we also write which of the above statements we use to prove  the existence of a $(x,y;u_1,u_2)$-set of chains in $B$. Later we give detailed arguments.

\begin{itemize}
\item[(a)] $k_1<k_2<k$, we use (i) for $i\in [k-1]\setminus \{k_1,k_2\}$,  (iii) for $j\in [2]$, and (ii). 
\item[(b)]  $k_1<k_2=k$, we use (i) for $i\in [k-1]\setminus \{k_1\}$,  (iii) for $j=1$, and (iv) for $j=2$.
\item[(c)]  $k_1<k<k_2$, we use (i) for $i\in [k-1]\setminus \{k_1\}$,   (iii) for $j=1$, and (ii). 
\item[(d)]  $k_1=k_2<k$, we use (i) for $i\in [k-1]\setminus \{k_1\}$,   (vi) for $j=1$, and (ii).   
\item[(e)] $k_1=k_2=k$, we use (i) for $i\in [k-1]$, and  (vii).    
\item[(f)] $k_1<k$ and $u_2=x$, we use  (i) for $i\in [k-1]\setminus \{k_1\}$, and (ii).
\item[(g)] $k_1=k$ and $u_2=x$,  we use  (i) for $i\in [k-1]$, and (v) for $j=1$. 
\end{itemize}

We prove cases (a), (c) and (g) in detail. Cases (b),(d) and (e) are similar to case (a), and case (f) is similar to case (g), so here we skip  details. 

{\em Case (a).} 
Suppose that $k_1<k_2<k$. By (i), there is a $(b_{i-1},b_i;b_{i-1},b_i)$-set of chains  ${\mathcal C_i}$ in 
$B_i$, for
$i\in [k-1]\setminus \{k_1,k_2\}$. By (iii), there is a $(b_{k_j-1},b_{k_j};b_{k_j-1},u_j)$-set of chains  ${\mathcal C_{k_j}}$  in $B_{k_j}$ for $j=1,2$. 
By (ii), there is  a $(b_{k-1},y;b_{k-1},b_k)$-set of chains  ${\mathcal C_k}$ in $B_k$.
Denote  by $P_{i}$ a ${\mathcal C_i}$-path in $B_{i}$, for $i\in [k]$.    

By Lemma \ref{skupek}, $G_0=B- \bigcup_{i=1}^{k}V(B_i)$ is a good chain  with respect to $x$.  Let $P_0$ be the path $x,b_0$. 
Define $P=\bigcup_{i=0}^{k} P_i$ (and recall that $P_i=B_i$, if $B_i$ is trivial). Then 
${\mathcal C}=\{G_0\}\cup \bigcup_{i=1}^k{\mathcal C_i}$ is a $(x,y;u_1,u_2)$-set of chains in $B$.   

{\em Case (c).} 
Suppose that $k_1<k<k_2$. 
 By (i) there is a $(b_{i-1},b_i;b_{i-1},b_i)$-set of chains  ${\mathcal C_i}$ in $B_i$, for
$i\in [k-1]\setminus \{k_1\}$. By (iii) there is a $(b_{k_1-1}b_{k_1};b_{k_1-1},u_1)$-set of chains ${\mathcal C_{k_1}}$ in $B_{k_1}$. 
By (ii) there is  a $(b_{k-1},y;b_{k-1},b_k)$-set of chains  ${\mathcal C_k}$ in $B_k$.
 
Since ${\mathcal C}_{k}$ is a $(b_{k-1},y;b_{k-1},b_k)$-set of chains in $B_k$, by Lemma \ref{spajanje} there is a  $(b_{k-1},y;b_{k-1},u_2)$-set of chains ${\mathcal D}_{k}$ in   $\bigcup_{i=k}^{k_2} B_i$.

By Lemma \ref{skupek}, $G_1=B- \bigcup_{i=1}^{k_2}V(B_i)$ is a good chain  with respect to $x$.  
Then $G_1$ together with chains in ${\mathcal C_i},i\in [k-1]$ and ${\mathcal D}_{k}$ forms a 
$(x,y;u_1,u_2)$-set of chains in $B$ (the corresponding path is  $P=\bigcup_{i=0}^{k} P_i$).

{\em Case (g).}  By (i) there is a $(b_{i-1},b_i;b_{i-1},b_i)$-set of chains  ${\mathcal C_i}$ in $B_i$, for
$i\in [k-1]$. Since $u_1\notin V(Q')$, by an assumption, we have $u_1\neq b_k$.  By (v), 
there is a $(b_{k-1},y;u_1,b_{k})$-set of chains  ${\mathcal C_k}$ in $B_k$.  By Lemma \ref{spajanje} there is a  $(b_{k-1},y;u_1)$-set of chains   ${\mathcal F}_{k}$ in $ \bigcup_{i=k}^{n} B_i$. 
   
Let $P_{i}$ be a ${\mathcal C_i}$-path in $B_{i}$, for $i\in [k]$, and define $P=\bigcup_{i=0}^{k} P_i$. Then chains in ${\mathcal C_i},i\in [k-1]$ and  ${\mathcal F}_{k}$ form a $(x,y;u_1,x)$-set of chains in $B$, with $P$ being the corresponding path.   
 \qed

\begin{definition}\label{defcik}
Let $B$ be a circuit graph with outer cycle $C$, and let   $u_1,u_2, u_3\in V(C)$. 
A set of pairwise disjoint chains ${\mathcal C}=\{G_1,\ldots,G_k\}$  is a $[u_1,u_2,u_3]$-set of chains in $B$ if there exists an 
even cycle $C'$ in $B$ such that 
\begin{itemize}
\item[(i)] $V(B)\setminus V(C')\subseteq \bigcup_{i=1}^{k}V(G_i)$ 
\item[(ii)] For $i\in [k]$, $G_i$ intersects $C'$ in exactly one vertex $x_i$, and 
$G_i$ is a good chain with respect to $x_i$. 
\item[(iii)]  For $j\in [3]$, either $G_i$ is a good chain with respect to $u_j$ and $x_i$ for some $i\in [k]$, or $u_j \notin \bigcup_{i=1}^{k}V(G_i)$.  
\end{itemize}
 A  cycle $C'$ that fulfills (i),(ii) and (iii)  is called a  ${\mathcal C}$-cycle. 
\end{definition}

If ${\mathcal C}$ fulfills (i),(ii) and (iii) for $j=1,2$, then ${\mathcal C}$ is a 
$[u_1,u_2]$-set of chains in $B$. 

\begin{lemma} \label{cikli}
Let $B$ be a bipartite circuit graph with outer cycle $C$, and let $u_1,u_2,u_3$ be any vertices of $C$. 
If all internal vertices of $B$ are of degree at least 4, then there exits a 
$[u_1,u_2,u_3]$-set of chains in $B$. 
\end{lemma}

\proof By Lemma \ref{plainchain}, $B-u_3$ is a plain chain of blocks 
$$B-u_3=B_1,b_1,B_2,\ldots,b_{n-1}, B_n\,.$$ 
Let $k_i\in [n]$ be such that $u_i\in V(B_i)\setminus V(B_{i-1})$ (here we set $B_0=\emptyset$).
 For $i\in [n]$,   define $P_i=B_i$ and ${\mathcal C_i}=\emptyset$, if $B_i$ is trivial. 
In the sequal we define $P_i$ and ${\mathcal C_i}$ for nontrivial blocks $B_i$.  

{\em Case 1: $k_1\neq k_2$}. By Lemma \ref{posebna} and \ref{to} there is 
  \begin{itemize}
\item[(i)] a $(b_{i-1},b_i;b_{i-1})$-set of chains ${\mathcal C_i}$ in $B_i$, if $i\notin \{k_1,k_2\}$, 
\item[(ii)] a $(b_{i-1},b_i;b_{i-1},u_j)$-set of chains ${\mathcal C_i}$ in $B_i$, if $i=k_j$ for $j=1,2$. 
\end{itemize}
Let   $P_i$ be the corresponding  ${\mathcal C_i}$-path in $B_i$, for $i\in [n]$. Let $b_0\in V(B_1)$  and $b_n\in V(B_n)$ be the neighbors of $x$ in $C$ and define $P_0=x,b_0$ and $P_{n+1}=b_n,x$. Define 
$C'=\bigcup_{i=0}^{n+1} P_i$ and 
$${\mathcal C}=\bigcup_{i=1}^n {\mathcal C_i}\,.$$
Then ${\mathcal C}$ is a $[u_1,u_2,u_3]$-set of chains, and $C'$ is a corresponding 
${\mathcal C}$-cycle.

{\em Case 2: $k_1= k_2$}. By Lemma \ref{posebna} and \ref{to} there is 
  \begin{itemize}
\item[(i)] a $(b_{i-1},b_i;b_{i-1},b_i)$-set of chains ${\mathcal C_i}$ in $B_i$, if $i\neq k_1$, 
\item[(ii)] a $(b_{k_1-1},b_{k_1};u_1,u_2)$-set of chains ${\mathcal C_{k_1}}$ in $B_{k_1}$. 
\end{itemize}
The rest of the proof is the same as above. 
\qed

\begin{lemma} \label{bipartite}
Let $B$ be a non-bipartite circuit graph with outer cycle $C$. Suppose that $x,y\in V(C)$ and that $Q$ is a $xy$-path in $C$. 
If all internal vertices of $B$ are of 
degree at least 4,  every vertex in $V(C)\setminus V(Q)$ is of degree at least 3 in $B$, and $B-x$ is bipartite, 
then for any vertex $u\in V(C-x)$  
there is an odd and an even $(x,y;u)$-set of chains in $B$.   
\end{lemma}

\proof
We claim that for any neighbor $z$ of $x$, there is a  
$(x,y;u)$-set of chains ${\mathcal C}$ in $B$ such that a ${\mathcal C}$-path contains the edge $xz$.    
 Before we prove the claim let us see how we prove the lemma using this claim. 
Since $B$ is non-bipartite and 2-connected, there is an odd cycle $C'$ containing $x$. Let $x_1$ and $x_2$ be the neigbors of $x$ in $C'$. Let $R$ be the $x_1x_2$-path in $C'$ not containing $x$.  
Suppose that $R_i$ is a $x_iy$-path in $B-x$ for $i=1,2$. Since $B-x$ is bipartite, $R_1\cup R_2\cup R$ is an even closed walk,  
and since $R$ is odd, $R_1$ and $R_2$ have different parities. 
It follows that every $x_1y$-path in $B-x$ is odd, and every $x_2y$-path in $B-x$ is even (or vice-versa). 
Using the above claim and setting $z=x_1$ (resp. $z=x_2$) we get an even (resp. an odd) $(x,y;u)$-set of chains in $B$.

In the rest of the proof we prove the claim. Let $$B-x=B_1,b_1,B_2,\ldots,b_{n-1}, B_n$$ 
and suppose that $xz\in E(B)$. By Lemma \ref{osnovna2}, $|V(B_i\cap Q)|\geq 2$ for $i\in [n]$,  
hence we may assume that $y\in V(B_n)\setminus V(B_{n-1})$. 
Let $k_u,k_{z}\in [n]$ be such that $u\in V(B_{k_u})\setminus V(B_{k_u+1})$ and  $z\in V(B_{k_{z}})\setminus V(B_{k_{z}+1})$ (here we set $B_{n+1}=\emptyset$). It follows from these definitions that $u\neq b_{k_u}$ and $z\neq b_{k_z}$. 
 
We shall construct  a $(x,y;u)$-set of chains ${\mathcal C}$ in $B$, and a ${\mathcal C}$-path $P$ in $B$,   so that $P$ contains the edge $xz$.  
We distinguish several cases with regard to $k_u$ and $k_z$. 
In each case we define $P_i=B_i$ and ${\mathcal C_i}=\emptyset$, if  $B_i$ is trivial and $i\in [n]$. Now we treat different cases and define 
$P_i$ and ${\mathcal C_i}$, if $B_i$ is nontrivial. 
 
 Suppose that $k_{z}<k_u<n$. By Lemma \ref{posebna} and Lemma \ref{to} there is 
\begin{itemize}
\item[(i)]  a $(z,b_{k_z};b_{{k_z}-1},b_{k_z})$-set of chains in $B_{k_z}$, where $k_z\neq k_u$
\item[(ii)] a $(b_{i-1},b_i;b_{i-1},b_i)$-set of chains in $B_i$, for $k_z< i< n,i\neq k_u$
\item[(iii)] a $(b_{n-1},y;b_{n-1})$-set of chains in $B_{n}$, where $n\neq k_u$
\item[(iv)] a $(b_{k_u-1},b_{k_u};u)$-set of chains in $B_{k_u}$.
\end{itemize}
 Let ${\mathcal C_i}$ be the set of chains in $B_i$  (as defined above), and let $P_i$ be  a  ${\mathcal C_i}$-path 
 for $i\in [n]\setminus [k_z-1]$.
By Lemma \ref{skupek}, $G_0=B-\bigcup_{i=k_z}^{n} V(B_i)$ is a good chain with respect to $x$ in $B$ (if 
$k_z=1$  this  is irrelevant).  
Let $P_0$ be the path $x,z$ and define $P=P_0\cup \bigcup_{i=k_z}^{n} P_i$. Then  
$${\mathcal C}=\{G_0\}  \cup \bigcup_{i=k_z}^{n} {\mathcal C_{i}}$$
is a $(x,y;u)$-set of chains in $B$ and $P$ is a ${\mathcal C}$-path.
If $k_z=k_u<n$ we use (ii) and (iii), and instead of (iv) we use 
\begin{itemize}
\item[(v)]  there is a $(z,b_{k_u};u)$-set of chains in $B_{k_u}$.
\end{itemize}
The rest of the proof is the same as above. 
If $k_z<k_u=n$ we use (i) and (ii), and instead of (iv) we use 
\begin{itemize}
\item[(vi)]  there is a $(b_{n-1},y;u)$-set of chains in $B_{n}$
\end{itemize}
and  the  rest of the proof is (again) the same as above (note that (v) and (vi) follow from  Lemma \ref{to}).

If $k_u<k_z<n$   then
 we use (i),(ii) and (iii). By Lemma \ref{spajanje} and (i),  
there is a $(z,b_{k_z}; b_{k_z},u)$-set of chains  ${\mathcal F_{k_z}}$ in 
$\bigcup_{i=k_u}^{k_z} B_i$. 
By Lemma \ref{skupek}, $G_1=B-\bigcup_{i=k_u}^{n} V(B_i)$ is a good chain with respect to $x$ in $B$. 
 It follows that 
 $${\mathcal C}=\{G_1\}\cup  {\mathcal F_{k_z}}  \cup \bigcup_{i=k_z+1}^{n} 
{\mathcal C_{i}}$$ 
 is a  
$(x,y;u)$-set of chains in $B$. The path $P$ (as defined above) is a ${\mathcal C}$-path.
This proves the claim when $k_z\neq n$.

Assume now that $k_z=n$.   
If $k_z=k_u=n$ and $z\neq y$ then, by   Lemma \ref{to}, there is
\begin{itemize}
\item[(vii)]  a $(z,y;u)$-set of chains  ${\mathcal H_{n}}$ in  $B_{n}$.    
\end{itemize}
 By Lemma \ref{skupek}, $G_{2}=B- V(B_n)$ 
is a good chain with respect to $x$ in $B$. It follows that
 ${\mathcal C}=\{G_{2}\}\cup {\mathcal H_{n}} $
 is a  
$(x,y;u)$-set of chains in $B$, and $P$ (as defined above) is a ${\mathcal C}$-path.

If $k_z=k_u=n$ and $z=y=u$ then $xz$ is an edge of $C$ (recall that $y\neq b_{n-1}$ and that 
$y$ is an external vertex of $B$). 
By Lemma \ref{osnovna}, $G_3=B-y$ is a good chain with respect to $x$.
Therefore 
${\mathcal C}=\{G_3\}$ is a $(x,y;u)$-set of chains,  where a ${\mathcal C}$-path in $B$ is the path  $x, y$. 
If $z=y\neq u$ then $B_n$ is a good chain with respect to $y$ and $u$, and $G_2$ is a good chain with respect to $x$. It follows that 
$\{B_n,G_2\}$ is a $(x,y;u)$-set of chains in $B$; again   a ${\mathcal C}$-path in $B$ is the path on two vertices $x,y$. If $z \neq y$ 

Finally, if $k_u<k_z=n$ and $z\neq y$ there is 
\begin{itemize}
\item[(viii)]  a $(z,y;b_{n-1})$-set of chains  in  $B_{n}$.    
\end{itemize}
 By Lemma \ref{spajanje} and (viii),  
there is a $(z,y;u)$-set of chains  ${\mathcal I_{n}}$ in 
$\bigcup_{i=k_u}^{n} B_i$. In this case 
 ${\mathcal C}=\{G_{1}\}\cup {\mathcal I_{n}}$
is a  $(x,y;u)$-set of chains in $B$. If $z=y$, then 
$G_4=\bigcup_{i=k_u}^{n} B_i$ is a good chain with respect to $u$ and $y$, hence 
 ${\mathcal C}=\{G_{1},G_4\} $
is a  $(x,y;u)$-set of chains in $B$. 
This proves the claim, and hence also the lemma. 
\qed

\begin{theorem}\label{main}
Let $B$ be a non-bipartite circuit graph with outer cycle $C$, and let  $x,y\in V(C)$. 
If all internal vertices of $B$ are of 
degree at least 4, 
then for any vertex $u\in V(C)$  
there is a  $(x,y;u)$-set of chains in $B$.   Moreover,  if  $Q$ is a $xy$-path in $C$ such that  
every vertex in $V(C)\setminus V(Q)$ is of degree at least 3 in $B$, then  for any vertex $u\in V(Q)$  
there is an odd and an even  $(x,y;u)$-set of chains in $B$. 
\end{theorem}

\proof  Suppose the theorem is not true. Let $B$ be a counterexample of minimum order. 
We may assume that $u\neq x$ (otherwise $u\neq y$, and the proof is analogous). 
By Lemma \ref{plainchain}, $B-x$ is a plain chain of blocks 
$$B-x=B_1,b_1,B_2,\ldots,b_{n-1}, B_n\,.$$
Let $k_1,k_2\in [n]$ be such that $u\in V(B_{k_1})\setminus V(B_{{k_1-1}})$ and $y\in V(B_{k_2})\setminus V(B_{k_2-1})$ 
(here we set $B_0=\emptyset$). 
  Let $b_0\in V(B_1)$  and $b_n\in V(B_n)$ be the neighbors of $x$ in $C$. We may assume that $xb_0$ is an edge of $Q$ and 
$xb_n$ is not an edge of $Q$, and 
that $u\in V(Q)$ (the last sentence of the theorem assumes $u\in V(Q)$, and for the proof of the 
first part of the theorem $u\in V(Q)$ may be assumed without loss of generality). Since $u\in V(Q)$ we have $k_1\leq k_2$.  
We give two constructions. In both constructions we define $P_i=B_i$ and  ${\mathcal C_{i}}=\emptyset$, 
if $B_i$ is trivial. 
In the sequal we treat nontrivial blocks $B_i$. \\

{\em Construction A.}  
If $k_1<k_2$ then, by minimality of $B$ (if $B_i$ is non-bipartite) and by Lemma \ref{to} (if $B_i$ is bipartite), there is 
\begin{itemize}
\item[(i)] a $(b_{i-1},b_i;b_i)$-set of chains in $B_i$, for $i\in [k_1-1]$,
\item[(ii)]  a $(b_{k_1-1},b_{k_1};u)$-set of chains  in $B_{k_1}$,
\item[(iii)]  a $(b_{i-1},b_i;b_{i-1})$-set of chains  in $B_i$,  for $i\in [k_2-1]\setminus [k_1]$,
\item[(iv)]  a $(b_{k_2-1},y;b_{k_2-1})$-set of chains in $B_{k_2}$,
\end{itemize}
  if $k_1=k_2$ and $y\neq b_0$ there is 
\begin{itemize}
 \item[(v)] a $(b_{k_2-1},y;u)$-set of chains in $B_{k_2}$.
\end{itemize}

Note that  for $i\in [k_2-1]$ every vertex in $V(C_i)\setminus V(Q_i)$ is of degree at least 3 in $B_i$, where $C_i$ 
is the outer cycle of $B_i$ and $Q_i=Q\cap B_i$. By minimality of $B$ we may apply the (last) statement of the theorem 
to $B_i$, if $B_i$ is non-bipartite. Hence, if $B_i$ is non-bipartite for some $i\in [k_2-1]$,   there is an odd and an even set of chains for (i), (ii) and (iii).  
Additionally, if $B_{k_2}$ is non-bipartite and $y=b_{k_2}$, then there is also an odd and an even set of chains
 ${\mathcal C_{k_2}}$ for  
(iv) and (v) (by minimality of $B$). 

Denote by ${\mathcal C_i}$ the set of chains in $B_i$ defined by (i)-(iv) if $k_1<k_2$; and defined  by 
(i) and (v) if $k_1=k_2$ and $y\neq b_0$. The ${\mathcal C_i}$-path is denoted by $P_i$, for $i\in [k_2]$. 
Let $P_0$ be the path $x,b_0$. 
Define $P=\bigcup_{i=0}^{k_2}P_i$. 
By Lemma \ref{skupek}, $G_1=B-\bigcup_{i=1}^{k_2}B_i$ is a good chain with respect to $x$ in $B$. 
Hence  
$${\mathcal C}=\{G_1\}\cup \bigcup_{i=1}^{k_2} {\mathcal C_i}$$
is a $(x,y;u)$-set of chains in $B$. The path $P$ is a ${\mathcal C}$-path in $B$. Moreover, if a block 
$B_i, i\in [k_2-1]$ is non-bipartite, then there exists an odd and an even  set of chains ${\mathcal C_i}$ in 
$B_i$, and so ${\mathcal C}$ is an odd or an even set of chains subject to the choice of  ${\mathcal C_i}$.  
  
If $y=b_0$ then $u=y=b_0$ (by our assumptions $u\in V(Q)$ and $u\neq x$). In this case 
$G_2=B-y$ is a good chain with respect to $x$, by Lemma \ref{osnovna}. Hence,  
${\mathcal C}=\{G_2\}$ is a $(x,y;u)$-set of chains, and $P=x,y$ is the corresponding ${\mathcal C}$-path. 
This proves the first claim of the theorem; and also the second claim of the theorem if $B_i$ 
is non-bipartite for some $i\in [k_2-1]$, or if $B_{k_2}$ is non-bipartite and $y=b_{k_2}$ 
(note that this, in particular, proves the theorem for the case if $y=b_n$ and $B-x$ is non-biparite). 
To finish the proof of the second claim of the theorem we give construction B, in which we assume that  
 every vertex in $V(C)\setminus V(Q)$ is of degree at least 3 in $B$.  We also assume that 
$B_i$ is  bipartite for  $i\in [k_2-1]$, and if  $y=b_{k_2}$  then $B_i$ 
is  bipartite for  $i\in [k_2]$.
\\

{\em Construction B.} If $k_1<k_2$ and $y\neq b_{k_2}$ then, by minimality of $B$ (if $B_i$ is non-bipartite) and by Lemma \ref{to} (if $B_i$ is bipartite),  there is 
\begin{itemize}
\item[(vi)]  a $(b_{i-1},b_i;b_{i-1})$-set of chains  in $B_i$, for $i\in [n]\setminus [k_2]$,
 \item[(vii)] a $(b_{k_2},y;b_{k_2-1})$-set of chains in $B_{k_2}$,
\end{itemize}
and if $k_1=k_2$  and $y\neq b_{k_2}$ there is 
 \begin{itemize}
\item[(viii)] a $(b_{k_2},y;u)$-set of chains in $B_{k_2}$.
\end{itemize}

Denote by ${\mathcal C_i}$ the set of chains in $B_i$ defined by (vi) and (vii) if $k_1<k_2$ and $y\neq b_{k_2}$;
and defined by (vi) and (viii) if $k_1=k_2$ and $y\neq b_{k_2}$. 
The ${\mathcal C_i}$-paths are denoted by $P_i$, for $i\in [n]\setminus [k_2-1]$. 

Suppose that $y\neq b_{k_2}$ and $k_1< k_2$. Let $R_1$ and $R_2$ be the $yb_{k_2}$-paths in $C_{k_2}$ (where 
$C_{k_2}$ is the outer cycle of $B_{k_2}$), and assume $b_{k_2-1}\in V(R_1)$. 
Since every vertex in $V(C)\setminus V(Q)$ is of degree at least 3 in 
$B$, we find that every vertex in 
$V(C_{k_2})\setminus V(R_1)$ is of degree at least 3 in $B_{k_2}$. Therefore, by minimailty of $B$, if $B_{k_2}$ is non-bipartite there is an odd and an  even set of chains 
${\mathcal C_{k_2}}$ for (vii) and (viii). Moreover, if $B_i$ is non-bipartite there exist  odd and  even sets of chains 
${\mathcal C_{i}}$ for $i\in [n]\setminus [k_2]$, as given by (vi). 

If 
$u\neq b_{k_1}$ then, by Lemma \ref{spajanje}  (see notes directly after Lemma 
\ref{spajanje}) and (vii), there is a $(b_{k_2},y;u)$-set of chains  ${\mathcal D_{k_2}}$  in $ \bigcup_{i=k_1}^{k_2} B_i$ (recall the assumption that $B_i$ is bipartite for  $i\in [k_2-1]$).
 By Lemma \ref{skupek}, $G_3=B-\bigcup_{i=k_1}^{n}B_i$ is a good chain with respect to $x$. 
Let $P_{n+1}$ be the path $x,b_n$. Define $P=\bigcup_{i=k_2}^{n+1}P_i$. 
 Then
$${\mathcal C}=\{G_3\}\cup  {\mathcal D_{k_2}}\cup \bigcup_{i=k_2+1}^{n} {\mathcal C_i}$$
is a $(x,y;u)$-set of chains in $B$.  The corresponding ${\mathcal C}$-path is $P$.

If  
$u=b_{k_1}$ the construction of a $(x,y;u)$-set of chains in $B$ is analogous as in the   case   
$u\neq b_{k_1}$ (the 
only difference is that ${\mathcal D_{k_2}}$ is a $(b_{k_2},y;u)$-set of chains in $ \bigcup_{i=k_1+1}^{k_2} B_i$, 
and  $G_3=B-\bigcup_{i=k_1+1}^{n}B_i$). 

If $y\neq b_{k_2}$ and  $k_1=k_2$, we use (vi) and (viii). In this case 
$${\mathcal C}=\{G_3\}\cup \bigcup_{i=k_2}^{n} {\mathcal C_i}$$ is a 
 $(x,y;u)$-set of chains in $B$.  
If $B_i$ is non-bipartite for some $i\in [n]\setminus [k_2-1]$ then we can choose ${\mathcal C}_i$ so that   ${\mathcal C}$ is an odd or an even  $(x,y;u)$-set of chains in $B$ (in all  cases above). 
This proves the second claim of the theorem if $y\neq b_{k_2}$ and $B_i$ is non-bipartite for some 
$i\in [n]\setminus [k_2-1]$. 

Suppose now that $y=b_{k_2},k_2\neq n$ and $u\neq b_{k_1}$.  If we use  (vi) for $i=k_2+1$, we find that   
${\mathcal C_{k_2+1}}$ is a $(b_{k_2},b_{k_2+1};b_{k_2})$-set of chains  in $B_{k_2+1}$. Hence, by   Lemma 
\ref{spajanje} (see   notes after Lemma 
\ref{spajanje}), there is a $(b_{k_2},b_{k_2+1};u)$-set of chains  ${\mathcal F_{k_2+1}}$ in 
 $ \bigcup_{i=k_1}^{k_2+1} B_i$ (recall the assumption that 
$B_i$ is  bipartite for  $i\in [k_2]$ if $y=b_{k_2}$). 
 Then
$${\mathcal C}=\{G_3\}\cup  {\mathcal F_{k_2+1}}\cup \bigcup_{i=k_2+2}^{n} {\mathcal C_i}$$  
is a $(x,y;u)$-set of chains in $B$. The corresponding ${\mathcal C}$-path is $P$.

If $y=b_{k_2}, k_2\neq n$ and $u=b_{k_1}$, then  let 
 ${\mathcal H_{k_2+1}}$ be a $(b_{k_2},b_{k_2+1};u)$-set of chains  in 
 $ \bigcup_{i=k_1+1}^{k_2+1} B_i$ (it exits by Lemma \ref{spajanje}) and 
define $G_4=B-\bigcup_{i=k_1+1}^{n} B_i$. 
 In this case 
 ${\mathcal C}=\{G_4\}\cup {\mathcal H_{k_2+1}}\cup    \bigcup_{i=k_2+2}^{n} {\mathcal C_i}$
is a $(x,y;u)$-set of chains in $B$.

Observe that, if  $B_i$  is non-bipartite for some $i\in  [n]\setminus[k_2]$, then  we can choose $P_i$, and hence also $P$,  so that 
${\mathcal C}$ is odd or even. This proves the second claim of the theorem if $y=b_{k_2}$ ($k_2\neq n$) and $B_i$ is non-bipartite for some  $i\in [n]\setminus [k_2]$. 

The last case to consider is when  $B_i$ is bipartite for $i\in [n]$. In this case $B-x$ is bipartite and the theorem follows from Lemma \ref{bipartite}. 
\qed

The proof of the following lemma  is similar to the proof of Lemma \ref{spajanje} (so we skip this proof).

\begin{lemma} \label{spajanje1}
 Let   $G$ be a bipartite plain chain of blocks 
$$G=B_1,b_1,B_2,\ldots,b_{n-1}, B_n\,.$$ 
Let $u\in V(B_j)\setminus \{b_{j-1},b_j\}$ for some $j\in \{2,\ldots,n-1\}$, and suppose that there exits a 
$[b_{j-1},b_j,u]$-set of chains 
in $B_j$.  
Then for any  $v'\in V(B_1)\setminus V(B_{2})$ and 
$v''\in V(B_n)\setminus V(B_{n-1})$, there exists a 
$[v',v'',u]$-set of chains in $G$. 
\end{lemma}

\begin{theorem} \label{glavni}
Let $B$ be a non-bipartite circuit graph with outer cycle $C$. Suppose that $x,y\in V(C)$ and that $B$ is good 
with respect to $x$ and $y$. If every internal vertex of $B$ is of degree at leat 4 and $B$ is not an odd cycle, 
 then there exists a $[x,y]$-set of chains in $B$. 
\end{theorem}

\proof 
By Lemma \ref{plainchain}, $B-x$ is a plain chain of blocks 
$$B-x=B_1,b_1,B_2,\ldots,b_{n-1}, B_n\,.$$ 
Let $b_0\in V(B_1)$  and $b_n\in V(B_n)$ be the neighbors of $x$ in $C$. 
Let $k \in [n]$ be such that $y\in V(B_{k})\setminus V(B_{{k-1}})$ (here we set $B_0=\emptyset$).

Suppose that at least one block $B_i$ is non-bipartite. If $B_i$ is trivial, define 
 ${\mathcal C_i}=\emptyset$ and $P_i=B_i$.  
Otherwise, by Theorem \ref{main} and Lemma \ref{to} there is  
\begin{itemize}
\item[(i)]  a $(b_{i-1},b_i;b_{i})$-set of chains  in $B_i$, for $i\in   [k-1]$,
 \item[(ii)] a $(b_{k-1},b_k;y)$-set of chains in $B_{k}$,
\item[(iii)] a $(b_{i-1},b_i;b_{i-1})$-set of chains in $B_{i}$ for $i\in [n]\setminus [k]$.
\end{itemize}

Denote by ${\mathcal C_i}$ the set of chains in $B_i$ defined by (i), (ii) and (iii). 
The ${\mathcal C_i}$-paths are denoted by $P_i$, for $i\in [n]$.
Let $P_0=x,b_0$ and $P_{n+1}=b_n,x$, and define $C'=\bigcup_{i=0}^{n+1} P_i$. 
Since  $B_i$ is non-bipartite for some $i\in [n]$,  there is an odd and an even set of chains 
  ${\mathcal C_i}$ for (i),(ii) or (iii). Hence, we can choose ${\mathcal C_i}$ and $P_i$ so that 
 $C'$ is even, and therefore ${\mathcal C}= \bigcup_{i=1}^n{\mathcal C_i}$ is a $[x,y]$-set of chains in $B$. 

Suppose that all blocks $B_i, i\in [n]$ are bipartite. 
Then all odd faces of $B$ are incident to $x$. Define $B_{n+1}=P_{n+1}$. 


If $y\notin\{b_{k-1},b_k\}$, then by Lemma \ref{cikli} there exits a  
$[b_{k-1},b_k,y]$-set of chains ${\mathcal D_k}$  in $B_k$. 
Let $G=B_1,b_1,B_2,\ldots,b_{n-1}, B_n,b_n,B_{n+1}\,.$ 
By Lemma \ref{spajanje1} there is a  $[x,y]$-set of chains in $G$, which is also  a 
$[x,y]$-set of chains in $B$ (because $G$ is a spanning subgraph of $B$).

Assume now that   $y\in\{b_{k-1},b_k\}$. 
Suppose that $y\in \{b_0,b_n\}$. We may assume $y=b_0$. 
If a block $B_i$ of $G$ is nontrivial, then 
by Lemma \ref{cikli}, there is a $[b_{i-1},b_i]$-set of chains ${\mathcal F_i}$ in $B_i$. Therefore, 
by Lemma \ref{spajanje1} there is a $[x,y]$-set of chains in $G$, which is also  a 
$[x,y]$-set of chains in $B$. Otherwise all blocks $B_i,i\in [n+1]$ are trivial. 
If $C$ is an even cycle, then $C$ itself is a $[x,y]$-set of chains in $B$. 
Otherwise $C$ is odd, and since $B$ is not an odd cycle,  $C$ has a chord. Hence $B$ has an even 
cycle  $C_0$ (which goes through $x$). Clearly, $C_0$ together with blocks $B_i$ such that 
$|V(B_i)\cap V(C_0)|\leq 1$   forms a $[x,y]$-set of chains in $B$.

Hence we may assume  that    $y=b_k$ where $k\notin\{0,n\}$. 
Suppose that all bounded odd faces of $B$ are incident to 
$y$ (and recall that all bounded odd faces are incident to $x$). 
Then there are exactly one or two such faces. However, if there is exactly one bounded odd face in $B$,  
and this odd face is incident to $x$ and $y$, then 
$xy\in E(C)$ (follows from the fact that $B$ is good with respect to $x$ and $y$) 
and so $y\in \{b_0,b_n\}$.

Therefore there are exactly  two bounded odd faces  in $B$ (both adjacent to   $x$ and $y$).
In this case the cycle $C'$(defined above) bounds exactly two 
odd faces of $B$ and therefore $C'$ is even.  
Hence, ${\mathcal C}$ (defined above) is a 
$[x,y]$-set of chains in $B$. 

We may therefore assume that there is a bounded odd face $F$ of $B$, which 
is not incident to $y$, and that $y\notin \{b_0,b_n\}$. 
 Let $xx_1$ and $xx_2$ be  edges incident to $F$.  Since $F$ is not incident to $y=b_k$ 
we may assume, without loss of generality, that $x_1,x_2\in \bigcup_{i=1}^{k} V(B_i)$.  
Since $F$ is an odd face and $G$ is bipartite, every $x_1x$-path in $G$ is odd and every 
$x_2x$-path in $G$ is even (or vice-versa). 
 
Let $k'\in [k]$ be such that $x_1\in V(B_{k'})\setminus V(B_{k'+1})$. 
If $B_{k'}$ resp. $B_i$ is nontrivial, then by Lemma \ref{posebna} and Lemma \ref{to} 
there is a 
\begin{itemize}
\item[(iv)]  a $(x_1,b_{k'};b_{k'-1},b_{k'})$-set of chains ${\mathcal G_{k'}}$ in $B_{k'}$,
\item[(v)] a $(b_{i-1},b_i;b_{i-1},b_i)$-set of chains ${\mathcal G_{i}}$ in $B_{i}$ for $i\in [n]\setminus [k']$.
\end{itemize}
 Let $P_i$ be the ${\mathcal G_{i}}$-path in $B_i$ for $i\in [n]\setminus [k'-1]$ (if $B_i$ is trivial, define $P_i=B_i$ and ${\mathcal G_{i}}=\emptyset$), and define 
$C''=\bigcup_{i=k'}^{n+1}P_i \cup \{xx_1\}$ 
(recall that $P_{n+1}=b_n,x$). Since every $x_1x$-path in $G$ is odd,  
$C''$ is even. 
By (iv) and Lemma \ref{spajanje}, there is a 
$(x_1,b_{k'};b_{k'})$-set of chains ${\mathcal H_{k'}}$ in $\bigcup_{i=1}^{k'}B_i$. Then  
${\mathcal G}=\bigcup_{i=k'+1}^n {\mathcal G_{i}}\cup {\mathcal H_{k'}}$
is a $[x,y]$-set of chains in $B$, and $C''$ is a ${\mathcal G}$-cycle in $B$. 
\qed

\begin{theorem}\label{final}
Let $B$ be a circuit graph such that every internal vertex of $B$ is of degree at least 4. 
Then $B$ has a spanning bipartite cactus. 
\end{theorem}
\proof
Let $C$ be the outer cycle of $B$. 
We prove a slightly stronger statement: if $B$ is a circuit graph such that every internal vertex 
of $B$ is of degree at least 4, and $x,y\in V(C)$ are vertices such that $B$ is good with respect to 
$x$ and $y$, then $B$ has a spanning bipartite cactus $T$ 
such that $x$ and $y$ are contained in exactly one block of $T$. 
  
The proof is by induction on $|V(B)|$.  The statement is clealy true if $B$ is an even cycle.  
If $B$ is an odd cycle and $B$ is good with respect to $x$ and $y$ then $x$ and $y$ are adjacent. A spanning $xy$-path in $B$ is a bipartite spanning cactus in $B$ such that $x$ and $y$ are contained in exactly one block of this cactus. 

If $B$ is not an odd cycle, then by Theorem \ref{glavni} (if $B$ is non-biparitite) and Lemma \ref{cikli} 
(if $B$ is bipartite), there is a 
$[x,y]$-set of chains ${\mathcal C}=\{G_1,\ldots,G_k\}$ in $B$. Let $C'$ be a 
${\mathcal C}$-cycle. 
Note that each block $B'$ of a chain $G_i,i\in [k]$ is good with respect to 
(both) cutvertices of $G_i$ contained in $B'$. 
Moreover, either $x\notin \bigcup_{i=1}^k V(G_i)$ or a chain of ${\mathcal C}$ 
is good with respect to $x$ (a similar fact is true for $y$).  
Therefore we can use the induction hypothesis, 
to obtain a spanning bipartite cactus $T(B')$ in $B'$ such that (both) cutvertices of  $G_i$ contained 
in $B'$
are contained in exactly one block of $T(B')$. Moreover the block $B_x$ that contains $x$ (if any) 
has a spanning bipartite cactus $T(B_x)$ such that $x$ is contained in exactly one block of $T(B_x)$        (a similar fact is true for $y$). Let ${\mathcal B}$ be the set of all blocks of  $G_i,i\in [k]$ and 
define  $T=C'\cup \bigcup_{B' \in {\mathcal B}}T(B')$. This gives the required bipartite cactus in $B$. 
\qed

Let ${\mathcal P}$ be the class of   3-connected planar graphs whose prisms are not hamiltonian.
We end the article with the problem to determine the minimum ratio of cubic vertices in a graph  
$G\in {\mathcal P}$.  Let $V_3(G)$ denote the set of cubic vertices in    $G$. 

\begin{problem}
Determine minimum $\epsilon$ such that there exist arbitrary large graphs $G\in {\mathcal P}$ with 
$|V_3(G)|/|V(G)|<\epsilon$. In particular, can $\epsilon$ be arbitrary small ?
\end{problem}

\noindent {\bf Acknowledgement:} 
The author thanks Uro\v s Milutinovi\'c for proofreading  parts of the final version of this paper. 
This work was supported by the Ministry of Education of Slovenia  [grant numbers P1-0297, J1-9109].

\end{document}